\documentclass[12pt]{article}
\usepackage{amssymb,amsmath,latexsym,geometry,cite,enumerate,float}
\newtheorem{theorem}{Theorem}
\newtheorem{lemma}[theorem]{Lemma}
\newtheorem{corollary}[theorem]{Corollary}

\usepackage{lineno}
\usepackage{setspace}

\begin{document}

%\linenumbers
\onehalfspace

\title{Some Comments on the Slater number}

\author{Michael Gentner and Dieter Rautenbach}

\date{}

\maketitle

\begin{center}
Institut f\"{u}r Optimierung und Operations Research, 
Universit\"{a}t Ulm, Ulm, Germany,
\{\texttt{michael.gentner, dieter.rautenbach}\}\texttt{@uni-ulm.de}\\[3mm]
\end{center}

\begin{abstract}
Let $G$ be a graph with degree sequence $d_1\geq \ldots \geq d_n$. Slater proposed $s\ell(G)=\min\{ s: (d_1+1)+\cdots+(d_s+1)\geq n\}$ as a lower bound on the domination number $\gamma(G)$ of $G$.
We show that deciding the equality of $\gamma(G)$ and $s\ell(G)$ 
for a given graph $G$ is NP-complete but that one can decide efficiently 
whether $\gamma(G)>s\ell(G)$ or $\gamma(G)\leq 
\left(\left\lceil\ln \left(\frac{n(G)}{s\ell(G)}\right)\right\rceil+1\right)s\ell(G)$.
For real numbers $\alpha$ and $\beta$ with $\alpha\geq \max\{ 0,\beta\}$,
let ${\cal G}(\alpha,\beta)$ be the class of non-null graphs $G$ 
such that every non-null subgraph $H$ of $G$ has at most $\alpha n(H)-\beta$ many edges.
Generalizing a result of Desormeaux, Haynes, and Henning,
we show that
$\gamma(G)\leq (2\alpha+1)s\ell(G)-2\beta$
for every graph $G$ in ${\cal G}(\alpha,\beta)$
with $\alpha \leq \frac{3}{2}$.
Furthermore, we show that 
$\gamma(G)/s\ell(G)$ 
is bounded for graphs $G$ in ${\cal G}(\alpha,\beta)$
if and only if $\alpha<2$.
For an outerplanar graph $G$ with $s\ell(G)\geq 2$, we show $\gamma(G)\leq 6s\ell(G)-6$.
In analogy to $s\ell(G)$, we propose 
$s\ell_t(G)=\min\{ s: d_1+\cdots+d_s\geq n\}$
as a lower bound on the total domination number.
Strengthening results due to Raczek
as well as Chellali and Haynes,
we show that
$s\ell_t(T)\geq \frac{n+2-n_1}{2}$
for every tree $T$ of order $n$ at least $2$ with $n_1$ endvertices.
\end{abstract}

{\small 

\begin{tabular}{lp{13cm}}
{\bf Keywords:} & 
Domination; 
Slater number; 
sparse graphs;
outerplanar graphs;
paired domination; 
total domination\\
{\bf MSC 2010:} & 05C69; 05C07
\end{tabular}
}

\pagebreak

\section{Introduction}

We consider finite, simple, and undirected graphs and use standard terminology.

One of the most well studied notion in graph theory is domination in graphs \cite{hhs}.
A set $D$ of vertices of a graph $G$ is a {\it dominating set} of $G$ if every vertex in $V(G)\setminus D$ has a neighbor in $D$,
where $V(G)$ is the vertex set of $G$.
The {\it domination number $\gamma(G)$} of $G$ is the minimum order of a dominating set of $G$.
Since the domination number is an NP-hard minimization parameter, 
upper bounds received more attention than lower bounds.
Slater \cite{s} proposed the following very simple lower bound merely depending on the degree sequence:
Let $G$ be a graph of order $n$ at least $1$, and let $d_1\geq \ldots \geq d_n$ be the non-increasing degree sequence of $G$.
The {\it Slater number $s\ell(G)$ of $G$} is the minimum positive integer $s$ for which $(d_1+1)+\cdots+(d_s+1)$ is at least $n$,
that is,
$$s\ell(G)=\min\Big\{ s: (d_1+1)+\cdots+(d_s+1)\geq n\Big\}.$$
Since the closed neighborhood $N_G[u]$ of a vertex $u$ of degree $d_G(u)$ contains exactly $d_G(u)+1$ elements,
the term $(d_1+1)+\cdots+(d_s+1)$ is an upper bound on the order of the union of the closed neighborhoods of any $s$ vertices of $G$,
which immediately implies 
\begin{eqnarray}\label{e1}
\gamma(G) & \geq & s\ell(G).
\end{eqnarray}
We present several results concerning algorithmic aspects of the Slater number,
its relation to the domination number for sparse graphs, 
and a variation of the Slater number for other domination parameters.
The next section contains our contributions 
together with a discussion of related results and references.

\section{Results}

Since (\ref{e1}) relies on a very simple argument, 
one might hope that the extremal graphs for this inequality have a simple structure.
Our first result shows that this hope is in vain.

\begin{theorem}\label{theorem1}
It is NP-complete to decide whether $\gamma(G)=s\ell(G)$ for a given graph $G$.
\end{theorem}
{\it Proof:} The proof relies on a reduction from {\sc 3-Sat} 
restricted to instances where every variable appears in at most five clauses (cf. [LO2] in \cite{gj}).
Therefore, let ${\cal F}$ be such an instance of {\sc 3-Sat} 
consisting of the clauses $c_1,\ldots,c_q$ 
over the boolean variables $x_1,\ldots,x_p$.
We construct a graph $G$ 
whose order is polynomially bounded in terms of $p$ and $q$ 
such that ${\cal F}$ is satisfiable if and only if $\gamma(G)=s\ell(G)$.

For every variable $x_i$, we create a clique $X_i$ of order $5p$, and select two special vertices $x_i$ and $\bar{x}_i$ within $X_i$.
Let $X=\{ x_1,\ldots,x_p,\bar{x}_1,\ldots,\bar{x}_p\}$.
For every clause $c_j$, we create a vertex $c_j$.
For every clause $c_j$ and every literal $y$ within the clause $c_j$,
we create an edge between the vertex $c_j$ and the vertex $y$.
This completes the construction of $G$.
Note that $G$ has order $n=5p^2+q$,
and that the vertices in $X$ are the $2p$ vertices of $G$ of largest degrees.
Let $d_1\geq \ldots \geq d_n$ be the degree sequence of $G$.

Since there are exactly $3q$ edges between $X$ and $\{ c_1,\ldots,c_q\}$,
the average degree of the vertices in $X$ is $5p-1+\frac{3q}{2p}$,
which implies
$(d_1+1)+\cdots+(d_p+1)\geq p\left( 5p+\frac{3q}{2p}\right)=5p^2+\frac{3}{2}q>n.$
Since every variable appears in at most five clauses,
the maximum degree of $G$ is at most $5p-1+5$,
which implies 
$(d_1+1)+\cdots+(d_{p-1}+1)\leq (p-1)(5p+5)=5p^2-5<n.$
These two inequalities imply that $s\ell(G)=p$.

If ${\cal F}$ is satisfiable, then the $p$ vertices in $X$ corresponding to the $p$ true literals form a dominating set,
which implies $\gamma(G)\leq p$. Since $\gamma(G)\geq s\ell(G)=p$, we obtain $\gamma(G)=s\ell(G)$.
Conversely, if $\gamma(G)=s\ell(G)$, then $s\ell(G)=p$, and the structure of $G$ imply that a minimum dominating set $D$ of $G$
contains exactly one vertex from each clique $X_i$.
Clearly, we may assume that $D\subseteq X$.
Since every vertex in $\{ c_1,\ldots,c_q\}$ has a neighbor in $D$,
the elements of $D$ indicate a satisfying truth assignment for ${\cal F}$,
which completes the proof. $\Box$

\medskip

\noindent Trivially, 
for every class of graphs, for which the domination number can be determined efficiently,
also the equality of the domination number and the Slater number can be decided efficiently.
This comment motivates the question whether 
deciding equality in (\ref{e1}) is still hard for chordal graphs.

If $G$ is a graph with degree sequence $d_1\geq \ldots\geq d_n$ such that 
$(d_1+1)+\cdots+(d_{s\ell(G)}+1)\leq n+k$, 
and every degree appears at most $k$ times within the degree sequence of $G$
for some fixed constant $k$,
then there are at most $f(k)\cdot n^{O(k)}$ many sets $S$ of exactly $s\ell(G)$ vertices of $G$
with $\sum\limits_{u\in S}d_G(u)\geq n$.
Since all these sets can be generated in polynomial time, 
also in this case equality in (\ref{e1}) can be decided efficiently.

Our next result investigates to which degree 
we can at least efficiently compare the domination number to the Slater number. 
Its proof relies on a result of Nemhauser and Wolsey \cite{nw} 
concerning the maximization of submodular functions.

For a graph $G$ of order $n$,
the function 
$$f:2^{V(G)}\to \mathbb{N}_0:D\mapsto \left|\bigcup\limits_{u\in D}N_G[u]\right|$$
is non-decreasing and submodular.
Let $D_0=\emptyset$, and, for every positive integer $i$ at most $n$, 
let $D_i$ be the set $D_{i-1}\cup \{ u_i\}$,
where $u_i\in V(G)\setminus D_{i-1}$ is chosen such that 
$f(D_{i-1}\cup \{ u_i\})$ is as large as possible,
that is, 
the vertices $u_1,u_2,\ldots,u_n$ are ordered greedily such that  
$\left|N_G[u_1]\cup \cdots \cup N_G[u_i]\right|$ grows quickly.
Since the selection of $u_i$ only depends on $f(D_{i-1}\cup \{ u_i\})$, 
such an ordering can be found in polynomial time.

For every two positive integers $\ell$ and $k$ at most $n$,
Nemhauser and Wolsey \cite{nw} showed that 
\begin{eqnarray}\label{e2}
f(D_\ell)\geq \left(1-e^{-\frac{\ell}{k}}\right)\max\Big\{ f(X):X\subseteq V(G)\mbox{ and }|X|=k\Big\}.
\end{eqnarray}

\begin{theorem}\label{theorem2}
It is possible to decide in polynomial time for a given graph $G$ whether 
$$\gamma(G)>s\ell(G)\mbox{ or }\gamma(G)\leq 
\left(\left\lceil\ln \left(\frac{n(G)}{s\ell(G)}\right)\right\rceil+1\right)s\ell(G).$$
\end{theorem}
{\it Proof:} Let $G$ be a graph of order $n$.
Let $s=s\ell(G)$, $p=\left\lceil\ln \left(\frac{n}{s}\right)\right\rceil$, and $\ell=ps$.
If $n<\ell$, then, trivially, $\gamma(G)\leq n<ps<\left(\left\lceil\ln \left(\frac{n}{s}\right)\right\rceil+1\right)s$.
Hence, we may assume that $n\geq \ell$.
Let $D_\ell$ be constructed greedily as above.
If $f(D_\ell)< \left(1-e^{-p}\right)n$, then (\ref{e2}) implies
$\max\Big\{ f(X):X\subseteq V(G)\mbox{ and }|X|=s\Big\}<n$,
which implies $\gamma(G)>s$.
If $f(D_\ell)\geq \left(1-e^{-p}\right)n$, 
then 
$D_\ell\cup \left(V(G)\setminus \bigcup\limits_{u\in D_\ell}N_G[u]\right)$
is a dominating set of $G$ of order at most 
$\ell+(n-f(D_\ell))\leq ps+\frac{n}{e^p}\leq \left(\left\lceil\ln \left(\frac{n}{s}\right)\right\rceil+1\right)s$,
which completes the proof.
$\Box$

\medskip

\noindent Requiring equality in (\ref{e1}) not only for a graph itself but also for all its induced subgraphs leads to a simple class of graphs that can easily be recognized.

\begin{theorem}\label{theorem5}
A graph $G$ satisfies $\gamma(H)=s\ell(H)$ for every induced subgraph $H$ of $G$
if and only if $G$ is $\{ K_2\cup K_1\cup K_1,C_4\cup K_1\}$-free.
\end{theorem}
{\it Proof:} Since 
$\gamma(K_2\cup K_1\cup K_1)>s\ell(K_2\cup K_1\cup K_1)$
and 
$\gamma(C_4\cup K_1)>s\ell(C_4\cup K_1)$,
the necessity follows.
In order to show the sufficiency, 
we may assume, for a contradiction,
that $G$ is a $\{ K_2\cup K_1\cup K_1,C_4\cup K_1\}$-free graph with $\gamma(G)>s\ell(G)$.
Since $\gamma(G)=1$ if and only if $s\ell(G)=1$, this implies $\gamma(G)\geq 3$.
Furthermore, the graph $G$ has at least one edge, say $xy$.
Let $z$ be a vertex in $V(G)\setminus (N_G[x]\cup N_G[y])$.
Since $G$ is $K_2\cup K_1\cup K_1$-free, the set $\{ x,y,z\}$ is dominating.
Since $\gamma(G)\geq 3$, there are vertices $x'$ and $y'$ with 
$N_G(x')\cap \{ x,y,z\}=\{ x\}$ and $N_G(y')\cap \{ x,y,z\}=\{ y\}$.
If $x'$ and $y'$ are not adjacent, then $\{ x,x',y',z\}$ induces $K_2\cup K_1\cup K_1$,
and, 
if $x'$ and $y'$ are adjacent, then $\{ x,y,x',y',z\}$ induces $C_4\cup K_1$,
which is a contradiction. 
$\Box$

\medskip

\noindent If $n$ is a positive even integer, and $G=K_{\frac{n}{2}}\cup \bar{K}_{\frac{n}{2}}$, 
then $s\ell(G)=2$ and $\gamma(G)=\frac{n}{2}+1$,
that is, in general, there is no upper bound on the domination number in terms of the Slater number.
For non-null trees $T$ though, 
Desormeaux, Haynes, and Henning \cite{dhh2}
showed $\gamma(T)\leq 3s\ell(T)-2$.
In \cite{ghr}, we showed $\gamma(G)\leq 3s\ell(G)+2k-2$ 
for graphs $G$ that arise by adding $k$ edges to $T$.
Our next results generalize this for sufficiently sparse graphs.

For real numbers $\alpha$ and $\beta$ with $\alpha\geq 0$,
let ${\cal G}(\alpha,\beta)$ be the class of non-null graphs $G$ 
such that every non-null subgraph $H$ of $G$ has at most $\alpha n(H)-\beta$ many edges.
Note that ${\cal G}(\alpha,\beta)$ is empty for $\beta>\alpha$.

\begin{theorem}\label{theorem3}
If $\alpha$ and $\beta$ are real numbers with $\max\{ 0,\beta\}\leq \alpha\leq \frac{3}{2}$, then 
$$\gamma(G)\leq (2\alpha+1)s\ell(G)-2\beta$$
for every graph $G$ in ${\cal G}(\alpha,\beta)$.
\end{theorem}
{\it Proof:} Let $G$ be a graph of order $n$ in ${\cal G}(\alpha,\beta)$.
Let $s=s\ell(G)$, and let $S$ be a set of $s$ vertices of $G$ of largest degrees.
Let $\Gamma=\sum\limits_{u\in S}d_G(u)$.
By the choice of $s$ and $S$, we have $\Gamma\geq n-s$.

Let $G'$ arise from the subgraph of $G$ induced by $\bigcup\limits_{u\in S}N_G[u]$
by removing all edges that are not incident with a vertex in $S$.
Let $G'$ have order $n'$.
Let $V_1=\{ u\in V(G')\setminus S:d_{G'}(u)=1\}$, and let $n_1=|V_1|$.
Let $n_2=n'-s-n_1$,
that is, $n_2$ is the number of vertices in $V(G')\setminus S$ that are of degree at least $2$ within $G'$.

Since $G''=G'-V_1$ is a non-null subgraph of $G$ of order $n'-n_1$, we obtain that 
$$m(G')=m(G'')+n_1\leq \alpha(n'-n_1)-\beta+n_1=\alpha n'-(\alpha-1)n_1-\beta.$$
Since the degree sum $2m(G')$ of $G'$ is at least $\Gamma+n_1+2n_2$,
we obtain
$$(n-s)+n_1+2n_2\leq \Gamma+n_1+2n_2\leq 2m(G')\leq 2\alpha n'-(2\alpha-2)n_1-2\beta.$$
Since $\alpha\leq \frac{3}{2}$, this implies
\begin{eqnarray*}
2\alpha n' & \geq & n-s+(2\alpha-1)n_1+2n_2+2\beta\\
& \geq & n-s+(2\alpha-1)(n_1+n_2)+2\beta\\
& = & n-s+(2\alpha-1)(n'-s)+2\beta,
\end{eqnarray*}
and, hence,
$n'\geq n-2\alpha s+2\beta$.
This implies that 
$S\cup (V(G)\setminus V(G'))$ is a dominating set of $G$ 
of order at most $s+(n-n')\leq (2\alpha+1)s-2\beta$,
which completes the proof.
$\Box$

\medskip

\noindent If $G$ is a non-null forest, then 
$G\in {\cal G}(1,1)$,
and Theorem \ref{theorem3} implies $\gamma(G)\leq 3s\ell(G)-2$,
the bound of Desormeaux, Haynes, and Henning \cite{dhh2} mentioned above.
If $G$ is a non-null cactus, that is, no two cycles of $G$ share an edge,
then 
$G\in {\cal G}\left(\frac{3}{2},\frac{3}{2}\right)$,
and Theorem \ref{theorem3} implies $\gamma(G)\leq 4s\ell(G)-3$.

For $\frac{3}{2}<\alpha<2$, 
the following weaker version of Theorem \ref{theorem3} still holds.

\begin{theorem}\label{theorem6}
If $\alpha$ and $\beta$ are real numbers with $\max\{ 0,\beta\}\leq \alpha<2$, then 
$$\gamma(G)\leq \left(\frac{5+2|\beta|}{2-\alpha}+5+3|\beta|\right)s\ell(G)$$
for every graph $G$ in ${\cal G}(\alpha,\beta)$.
\end{theorem}
{\it Proof:} Let $G$ be a graph of order $n$ in ${\cal G}(\alpha,\beta)$.
Let $s=s\ell(G)$, and let $S$ be a set of $s$ vertices of $G$ of largest degrees.
Let $T=\left(\bigcup\limits_{u\in S}N_G[u]\right)\setminus S$, and let $R=V(G)\setminus (S\cup T)$.
Let $t=|T|$, and $r=|R|$.
Let $m_s$ be the number of edges of the subgraph of $G$ induced by $S$,
and let $m_t$ be the number of edges of $G$ between $S$ and $T$.

Let $\epsilon=2-\alpha$.

For a contradiction, suppose that 
$\gamma(G)>\left(\frac{5+2|\beta|}{2-\alpha}+5+3|\beta|\right)s\ell(G)
=\left(\frac{5+2|\beta|}{\epsilon}+5+3|\beta|\right)s$.
Since $S\cup R$ is a dominating set of $G$, we obtain $\gamma(G)\leq s+r$, 
and, hence,
$$r\geq \gamma(G)-s>\left(\frac{5+2|\beta|}{\epsilon}+4+3|\beta|\right)s.$$
Since $G\in {\cal G}(\alpha,\beta)$, we have $m_s\leq \alpha s-\beta\leq 2s+|\beta|$.

By the choice of $s$ and $S$, 
we have 
$$2m_s+m_t=\sum\limits_{u\in S}d_G(u)\geq n-s=t+r.$$
Since $s\geq 1$, this implies
\begin{eqnarray*}
m_t & \geq  & t+r-2m_s\\
&\geq & t+r-4s-2|\beta|\\
&>& t+\left(\frac{5+2|\beta|}{\epsilon} +4+3|\beta|\right)s-4s-2|\beta|\\
& \geq & t+\left(\frac{5+2|\beta|}{\epsilon}+|\beta|\right)s.
\end{eqnarray*}
Since $m_t$ is an integer, we obtain
$m_t\geq t+t'$,
where $t'=\left\lceil 
\left(\frac{5+2|\beta|}{\epsilon}+|\beta|\right)s 
\right\rceil$.

We consider three cases.

\medskip

\noindent {\bf Case 1} $t\geq t'$.

\medskip

\noindent There is a set $T'$ of $t'$ vertices in $T$ such that there are at least $2t'$ edges between $S$ and $T'$. 
This implies that the subgraph $H$ of $G$ induced by $S\cup T'$ satisfies
\begin{eqnarray*}
\frac{m(H)+\beta}{n(H)} 
& \geq & \frac{m(H)-|\beta|}{n(H)}\\
& \geq & \frac{2t'-|\beta|}{s+t'}\\
& \geq & \frac{2\left(\frac{5+2|\beta|}{\epsilon}+|\beta|\right)s-|\beta|}{s+\left(\frac{5+2|\beta|}{\epsilon}+|\beta|\right)s+1}\\
& \geq & \frac{\left(\frac{10+4|\beta|}{\epsilon}+|\beta|\right)s}{\left(\frac{5+2|\beta|}{\epsilon}+|\beta|+2\right)s}\\
& > & 2-\epsilon,
\end{eqnarray*}
which is a contradiction.

\medskip

\noindent {\bf Case 2} $t<t'$ and $t\geq \frac{(2-\epsilon)s+|\beta|}{\epsilon}$.

\medskip

\noindent Since there are more than $2t$ edges between $S$ and $T$,
the subgraph $H$ of $G$ induced by $S\cup T$ satisfies
\begin{eqnarray*}
\frac{m(H)+\beta}{n(H)} 
& \geq & \frac{m(H)-|\beta|}{n(H)}\\
& > & \frac{2t-|\beta|}{s+t}\\
& \geq & 2-\epsilon,
\end{eqnarray*}
which is a contradiction.

\medskip

\noindent {\bf Case 3} $t<t'$ and $t<\frac{(2-\epsilon)s+|\beta|}{\epsilon}$.

\medskip

\noindent Since there are more than $t'$ edges between $S$ and $T$,
the subgraph $H$ of $G$ induced by $S\cup T$ satisfies
\begin{eqnarray*}
\frac{m(H)+\beta}{n(H)} 
& \geq & \frac{m(H)-|\beta|}{n(H)}\\
& > & \frac{t'-|\beta|}{s+t}\\
& > & \frac{\left(\frac{5+2|\beta|}{\epsilon}+|\beta|\right)s-|\beta|}{s+\frac{(2-\epsilon)s+|\beta|}{\epsilon}}\\
& \geq & \frac{\left(\frac{5+2|\beta|}{\epsilon}\right)s}{\frac{2s+|\beta|}{\epsilon}}\\
& \geq & \frac{5s+2|\beta|s}{2s+|\beta|s}\\
& > & 2,
\end{eqnarray*}
which is a contradiction.

This completes the proof.
$\Box$

\medskip

\noindent For $\alpha\geq 2$, 
the fraction of the domination number and the Slater number
is no longer bounded within ${\cal G}(\alpha,\beta)$.
The fundamentally different behaviour 
for $\alpha<2$
and 
$\alpha\geq 2$
is reflected by the following result.

\begin{corollary}\label{corollary1}
If $\alpha$ and $\beta$ are real numbers with $\alpha\geq \max\{ 0,\beta\}$, then 
$\sup\left\{ \frac{\gamma(G)}{s\ell(G)}:G\in {\cal G}(\alpha,\beta)\right\}<\infty$
if and only if $\alpha<2$.
\end{corollary}
{\it Proof:} If $\alpha<2$, then Theorem \ref{theorem6} implies that the fraction $\gamma(G)/s\ell(G)$ is bounded for the graphs $G$ in ${\cal G}(\alpha,\beta)$.
For $\alpha\geq 2$, and positive even $n$, the graph $G^*=K_{2,\frac{n}{2}-1}\cup \bar{K}_{\frac{n}{2}-1}$ belongs to ${\cal G}(2,2)$,
and, hence, also to ${\cal G}(\alpha,\beta)$.
Since $s\ell(G^*)=2$ and $\gamma(G^*)=\frac{n}{2}+1$,
the fraction $\gamma(G)/s\ell(G)$ is unbounded within ${\cal G}(\alpha,\beta)$.
$\Box$

\medskip

\noindent Since the graph $G^*$ considered in the previous proof is planar,
the fraction $\gamma(G)/s\ell(G)$ is unbounded for planar graphs.
Since $G^*$ is not outerplanar though, this could be different for outerplanar graphs.
Note that outerplanar graphs of order $n$ at least $2$ may have up to $2n-3$ edges, 
which implies that Theorem \ref{theorem6} does not apply to them.
Therefore, in order to show the boundedness of $\gamma(G)/s\ell(G)$ for these graphs,
which is our next goal, we need to exploit more than their density.

\begin{lemma}\label{lemma1}
Let $G$ be an outerplanar graph,
and let $S$ be a set of $s$ vertices of $G$ with $s\geq 2$.
Let $T_2$ be the set of vertices in $V(G)\setminus S$ 
that have at least two neighbors in $S$.

If $t_2=|T_2|$,
$m_s$ is the number of edges of the subgraph of $G$ induced by $S$,
and $m_2$ is the number of edges between $S$ and $T_2$,
then $2m_s+m_2-t_2\leq 5s-6$.
\end{lemma}
{\it Proof:} Clearly, we may assume that $V(G)=S\cup T_2$, 
and that $G$ is maximal outerplanar.
The proof is by induction on the order $n$ of $G$.
For $n\leq 4$, the statement is easily verified.
Now, let $n\geq 5$.
Let $G$ be embedded in the plane such that all vertices lie on the boundary of the unbounded face,
and let $C$ be the Hamiltonian cycle of $G$ forming the boundary of that face.

First, we assume that some vertex $u$ in $T_2$ has a neighbor $v$ in $S$ 
such that the edge $uv$ does not belong to $C$.
The graph $G$ is the union of two maximal outerplanar graphs 
$G^{(1)}$ and $G^{(2)}$ of orders at least $3$ but less than $n$ 
such that $G^{(1)}$ and $G^{(2)}$ share exactly the edge $uv$.
For $i\in \{ 1,2\}$,
let $S^{(i)}=S\cap V(G^{(i)})$, 
let $s^{(i)}=|S^{(i)}|$,
and let $m_s^{(i)}$ be the number of edges of the subgraph of $G^{(i)}$ induced by $S^{(i)}$.
Since $S^{(1)}\cup S^{(2)}=S$ and $S^{(1)}\cap S^{(2)}=\{ v\}$, 
we obtain 
$s^{(1)}+s^{(2)}=s+1$
and
$m_s=m_s^{(1)}+m_s^{(2)}$.
If $s^{(i)}=1$, 
then the vertices in the non-empty set $V(G^{(i)})\setminus \{ u,v\}$ 
belong neither to $S$ nor to $T_2$, which contradicts $V(G)=S\cup T_2$.
Hence, $s^{(1)},s^{(2)}\geq 2$.
For $i\in \{ 1,2\}$,
let $T_2^{(i)}$ be the set of vertices in $V(G^{(i)})\setminus S^{(i)}$ 
that have at least two neighbors in $S^{(i)}$, 
let $t_2^{(i)}=|T_2^{(i)}|$,
and let $m_2^{(i)}$ be the number of edges between $S^{(i)}$ and $T_2^{(i)}$.
Note that $T_2^{(1)}\cup T_2^{(2)}=T_2$ and $T_2^{(1)}\cap T_2^{(2)}\subseteq \{ u\}$.
If $T_2^{(1)}\cap T_2^{(2)}=\emptyset$, then 
$t_2^{(1)}+t_2^{(2)}=t_2$
and 
$m_2^{(1)}+m_2^{(2)}=m_2$.
If $T_2^{(1)}\cap T_2^{(2)}=\{ u\}$, then 
$t_2^{(1)}+t_2^{(2)}=t_2+1$
and
$m_2^{(1)}+m_2^{(2)}=m_2+1$,
because the edge $uv$ contributes to $m_2^{(1)}$ as well as to $m_2^{(2)}$.
In both cases, 
$m_2\leq m_2^{(1)}+m_2^{(2)}$
and
$t_2\geq t_2^{(1)}+t_2^{(2)}-1$,
which, by induction, implies
\begin{eqnarray*}
2m_s+m_2-t_2 & \leq & 2(m_s^{(1)}+m_s^{(2)})+(m_2^{(1)}+m_2^{(2)})-(t_2^{(1)}+t_2^{(2)}-1)\\
&=& (2m_s^{(1)}+m_2^{(1)}-t_2^{(1)})+(2m_s^{(2)}+m_2^{(2)}-t_2^{(2)})+1\\
& \leq & (5s_1-6)+(5s_2-6)+1\\
& = & 5(s_1+s_2-1)-6\\
& = & 5s-6.
\end{eqnarray*}
Next, we may assume that all edges between $S$ and $T_2$ are edges of $C$.
This immediately implies $m_2-t_2\leq s$.
Furthermore, since $S$ induces an outerplanar graph of order at least $2$, we have $m_s\leq 2s-3$.
Altogether, we obtain $2m_s+m_2-t_2\leq 4s-6+s=5s-6$,
which completes the proof. $\Box$

\begin{theorem}\label{theorem7}
If $G$ is an outerplanar graph with $s\ell(G)\geq 2$, then $\gamma(G)\leq 6s\ell(G)-6$.
\end{theorem}
{\it Proof:} Let $n$, $S$, $s$, $T$, $t$, $R$, $r$, and $m_s$ be as in the proof of Theorem \ref{theorem6}.
Let $T_2$ be the set of vertices in $T$ that have at least two neighbors in $S$,
and let $m_2$ be the number of edges between $S$ and $T_2$.
Lemma \ref{lemma1} implies $2m_s+m_2\leq 5s-6+|T_2|$.
This implies that there are at least 
$\sum\limits_ {u\in S}d_G(u)-2m_s-m_2\geq (n-s)-(5s-6+|T_2|)=n-6s+6-|T_2|$ edges between $S$ and $T\setminus T_2$.
Since every vertex in $T\setminus T_2$ has exactly one neighbor in $S$,
this implies $s+t=|S|+|T_2|+|T\setminus T_2|\geq s+|T_2|+(n-6s+6-|T_2|)=n-5s+6$,
and, hence, $r=n-(s+t)\leq 5s-6$.
Since $S\cup R$ is a dominating set of $G$, we obtain $\gamma(G)\leq s+r\leq 6s-6$,
which completes the proof. $\Box$

\medskip

\noindent We believe that Theorem \ref{theorem7} can still be improved a little bit.
Let $s$ be an even integer at least $4$.
Let $G_S$ be a maximal outerplanar graph of order $s$ embedded such the Hamiltonian cycle
$u_1u_2\ldots u_su_1$ forms the boundary of the unbounded face.
If $G$ arises from $G_S$ by 
\begin{itemize}
\item adding $s$ vertices $v_1,\ldots,v_s$, 
where $v_i$ is adjacent to $u_i$ and $u_{i+1}$, and induces are taken modulo $s$,
and 
\item adding further $5s-6$ isolated vertices,
\end{itemize}
then $G$ is outerplanar,
$s\ell(G)=s$,
and 
$\gamma(G)$ is at least $\frac{s}{2}+(5s-6)$.
This example suggests that the factor ``$6$'' might be replaced by ``$11/2$'' but not by less.
Further improvement seems possible for maximal outerplanar graphs.

\medskip

\noindent For many variants of the domination number, lower bounds that are similar to the Slater number can be defined.
As an example we consider paired domination \cite{hs} and total domination \cite{hy}.

Let $G$ be a graph.
A set $D$ of vertices of $G$ is a {\it paired dominating set} of $G$ if $D$ is a dominating set, 
and the subgraph of $G$ induced by $D$ has a perfect matching.
The {\it paired domination number $\gamma_p(G)$} of $G$ is the minimum order of a paired dominating set of $G$.
A set $D$ of vertices of $G$ is a {\it total dominating set} of $G$ if every vertex of $G$ has a neighbor in $G$.
The {\it total domination number $\gamma_t(G)$} of $G$ is the minimum order of a total dominating set of $G$.

If $G$ has order $n$ at least $1$, and $d_1\geq \ldots \geq d_n$ is its non-increasing degree sequence,
then let
$$s\ell_t(G)=\min\Big\{ s: d_1+\cdots+d_s\geq n\Big\}.$$
Obviously,
$\gamma_p(G) \geq \gamma_t(G) \geq s\ell_t(G).$
Let $T$ be a tree of order $n$ at least $3$ with $n_1$ endvertices.
Raczek \cite{r} proved $\gamma_p(T) \geq \frac{n+2-n_1}{2}$,
and Chellali and Haynes \cite{ch} proved the stronger result $\gamma_t(T) \geq \frac{n+2-n_1}{2}$.
Both these results were inspired by Lema\'{n}ska's \cite{l} lower bound 
$\gamma(T)  \geq \frac{n+2-n_1}{3}$ on the domination number of $T$.
As observed by Desormeaux, Haynes, and Henning \cite{dhh2}, 
the right hand side of the previous inequality is actually a lower bound for the Slater number rather than the domination number, that is,
$s\ell(T) \geq \frac{n+2-n_1}{3}$, which strengthens Lema\'{n}ska's result.
Our final result strengthens the above bounds due to Raczek \cite{r} and Chellali and Haynes \cite{ch} in a similar way.

\begin{theorem}\label{theorem4}
If $T$ is a tree of order $n$ at least $2$ with $n_1$ endvertices, then $s\ell_t(T)\geq \frac{n+2-n_1}{2}$.
\end{theorem}
{\it Proof:} Let $d_1\geq \ldots \geq d_n$ be the non-increasing degree sequence of $T$.
Let $s=s\ell_t(T)$.
Let $V_1=\{ u\in V(T):d_T(u)=1\}$, and let $V_2=V(T)\setminus V_1$.
If $T$ is a star, then $s\ell_t(T)=2>\frac{3}{2}=\frac{n+2-n_1}{2}$.
Hence, we may assume that $T$ is not a star, which implies $n_1\leq n-2$.
Since $T$ has degree sum $n_1+\sum\limits_{u\in V_2}d_G(u)$, and exactly $n-1$ edges, 
we obtain $n_1+\sum\limits_{u\in V_2}d_G(u)=2n-2$.
This implies $\sum\limits_{u\in V_2}d_G(u)\geq 2n-2-n_1\geq n$,
and, hence, $d_s\geq 2$.
Since $T$ is a tree of order at least $2$, we have $n_1=\sum\limits_{u\in V_2}(d_G(u)-2)+2$.
Now, 
$n
\leq \sum\limits_{i=1}^s d_i
=\sum\limits_{i=1}^s (d_i-2)+2s
\leq \sum\limits_{u\in V_2}(d_G(u)-2)+2s
=n_1-2+2s$, 
which implies $s\geq \frac{n+2-n_1}{2}$. $\Box$

\medskip 

\noindent Similarly as in \cite{dhh2}, it is easy to see that equality holds in the above theorem 
if and only if $n$ and $n_1$ have the same parity modulo $2$, and $d_{s\ell_t(T)+1}\leq 2$.
Simple modifications of the proof of Theorem \ref{theorem1} 
imply that it is NP-complete to decide
whether $\gamma_t(G)=s\ell_t(G)$
for a given graph $G$.
Also the above results concerning sparse graphs can be extended to $s\ell_t(G)$
and the paired/total domination number.

In \cite{dhh1}, Desormeaux, Haynes, and Henning define the {\it connected order-sum number} ${\rm ord}_c(G)$
of a graph $G$ with non-increasing degree sequence $d_1\geq \ldots \geq d_n$ as 
${\rm ord}_c(G)=\min\Big\{s: d_1+\cdots+d_s\geq n-s+2\Big\}.$
They show that ${\rm ord}_c(G)$ is a lower bound on the connected domination number of $G$
and that equality holds for all trees.


\begin{thebibliography}{}
\bibitem{ch} M. Chellali, T.W. Haynes, A note on the total domination number of a tree, Journal of Combinatorial Mathematics and Combinatorial Computing 58 (2006) 189-193.
\bibitem{dhh1} W.J. Desormeaux, T.W. Haynes, M.A. Henning, Bounds on the connected domination number of a graph, Discrete Applied Mathematics 161 (2013) 2925-2931.
\bibitem{dhh2} W.J. Desormeaux, T.W. Haynes, M.A. Henning, Improved bounds on the domination number of a tree, Discrete Applied Mathematics 177 (2014) 88-94.
\bibitem{gj} M.R. Garey, D.S. Johnson, Computers and Intractability: A Guide to the Theory of NP-Completeness, W.H. Freeman \& Co. New York, New York, 1979.
\bibitem{ghr} M. Gentner, M.A. Henning, D. Rautenbach, Smallest Domination Number and Largest Independence Number of Graphs and Forests with given Degree Sequence, arXiv:1507.04647v1
\bibitem{hhs} T.W. Haynes, S.T. Hedetniemi, P.J. Slater, Fundamentals of Domination in Graphs, Marcel Dekker, Inc., New York, 1998.
\bibitem{hs} T.W. Haynes, P.J. Slater, Paired-domination and the paired-domatic number, Congressus Numerantium 109 (1995) 65-72.
\bibitem{hy} M.A. Henning, A. Yeo, Total domination in graphs, Springer, 2013.
\bibitem{l} M. Lema\'{n}ska, Lower bound on the domination number of a tree, Discussiones Mathematicae Graph Theory 24 (2004) 165-169.
\bibitem{nw} G.L. Nemhauser, L.A. Wolsey, Best Algorithms for Approximating the Maximum of a Submodular Set Function, Mathematics of Operations Research 3 (1978) 177-188.
\bibitem{r} J. Raczek, Lower bound on the paired domination number of a tree, Australasian Journal of Combinatorics 34 (2006) 343-347.
\bibitem{s} P.J. Slater, Locating dominating sets and locating-dominating sets, in: Graph Theory, Combinatorics, and Applications: Proc. 7th Quadrennial Int. Conf. Theory Applic. Graphs 2 (1995), 1073-1079.
\end{thebibliography}
\end{document}